\documentclass[10pt]{article}
\usepackage{epsfig,amssymb}
\usepackage{palatino}
\usepackage{verbatim}
\usepackage{amsmath}

\def\baa {\begin{eqnarray*}}
\def\eaa {\end{eqnarray*}}

\def \la {\lambda}

\def\Frac#1#2{\mbox{\large${\textstyle \frac{#1}{#2}}$}}

\def \ra {{\quad\Rightarrow\quad}}

\def \sign {{\rm sign\,}}
\def \Ai {{\rm Ai}}

\def\N{{\mathbb N}}

\def\OO{{\cal O}}

\def\la{\lambda}

\newtheorem{lemma}{Lemma}[section]
\newtheorem{proposition}[lemma]{Proposition}
\newtheorem{corollary}[lemma]{Corollary}
\newtheorem{theorem}[lemma]{Theorem}
\newtheorem{remark}[lemma]{Remark}

\newtheorem{definition}[lemma]{Definition}

\def\be  {\begin{equation}}
\def\ee  {\end{equation}}
\def\ba  {\begin{eqnarray}}
\def\ea  {\end{eqnarray}}
\def\baa {\begin{eqnarray*}}
\def\eaa {\end{eqnarray*}}
\def\bc  {}
\makeatletter
\@addtoreset{equation}{section}
\makeatother
\def\proof{\medskip\noindent{\bf Proof.} }
\def\qed{\hfill $\Box$}
\newcommand {\lb} {\label}
\newcommand {\rf[1]} {(\ref{#1})}

\begin{document}

\title{On the largest critical value of $T_n^{(k)}$ }

\author{N.\,Naidenov, G.\,Nikolov, A.\,Shadrin}

\date{October 16, 2017}
\maketitle



\begin{abstract}
We study the quantity
$$
\tau_{n,k}:=\frac{|T_n^{(k)}(\omega_{n,k})|}{T_n^{(k)}(1)}\,,
$$
where $T_n$ is the Chebyshev polynomial of degree $n$, and
$\omega_{n,k}$ is the rightmost zero of $T_n^{(k+1)}$.

Since the absolute values of the local maxima of $T_n^{(k)}$
increase monotonically towards the end-points of $[-1,1]$, the value
$\tau_{n,k}$ shows how small is the largest critical value of
$\,T_n^{(k)}\,$ relative to its global maximum $\,T_n^{(k)}(1)$.

In this paper, we improve and extend earlier estimates by
Erd\H{o}s--Szeg\H{o},
Eriksson and Nikolov in several directions.

Firstly, we show that the sequence $\,\{\tau_{n,k}\}_{n=k+2}^{\infty}$ is
monotonically decreasing in $n$, hence derive several sharp estimates,
in particular
$$
    \tau_{n,k} \le \left\{
\begin{array}{rcll}
   \tau_{k+4,k}
&=& \frac{1}{2k+1}\,\frac{3}{k+3}\,, & n \ge k+4\,, \\
   \tau_{k+6,k}
&=& \frac{1}{2k+1}\, (\frac{5}{k+5})^2 \beta_k\,,
    & n \ge k+6\,,
\end{array} \right.
$$
where $\beta_k < \frac{2+\sqrt{10}}{5} \approx 1.032$.

We also obtain an upper bound which is uniform in $n$ and $k$,
and that implies in particular
$$
   \tau_{n,k} \approx \big(\Frac{2}{e}\big)^k, \quad n \ge k^{3/2}; \qquad
   \tau_{n,n-m} \approx \big(\Frac{em}{2}\big)^{m/2} n^{-m/2}; \qquad
   \tau_{n,n/2} \approx \big(\Frac{4}{\sqrt{27}}\big)^{n/2}.
$$
Finally, we derive the exact asymptotic formulae for the quantities
$$
\tau_k^{*} := \lim_{n\to\infty}\tau_{n,k} \quad \mbox{ and }\quad
\tau_m^{**} := \lim_{n\to\infty} n^{m/2} \tau_{n,n-m}\,,
$$
which show that our upper bounds for $\tau_{n,k}$ and $\tau_{n,n-m}$
are asymptotically correct
with respect to the exponential terms given above.
\end{abstract}


\section{Introduction and statement of the results}


We study the quantity
$$
\tau_{n,k}:=\frac{|T_n^{(k)}(\omega_{k})|}{T_n^{(k)}(1)}\,,
$$
where $T_n$ is the Chebyshev polynomial of degree $n$, and
$\omega_k$ is the rightmost zero of $T_n^{(k+1)}$.

Since the absolute values of the local maxima of $T_n^{(k)}$
increase monotonically towards the end-points of $[-1,1]$, the value
$\tau_{n,k}$ shows how small is the largest critical value of
$\,T_n^{(k)}\,$ relative to its global maximum $\,T_n^{(k)}(1)$ (see
Figure~1).

This value is useful in several applications which include
some Markov-type inequalities \cite{es},
\cite{n1}, \cite{s2}, the Landau--Kolmogorov inequalities for
intermediate derivatives \cite{e}, \cite{s2}, where the estimates of
$f^{(k)}$ on a  subinterval slightly smaller than $[-1,1]$ are
needed, and also in studying extreme zeros of ultraspherical polynomials.

\begin{figure}[h]
\centering
\includegraphics[scale=1]{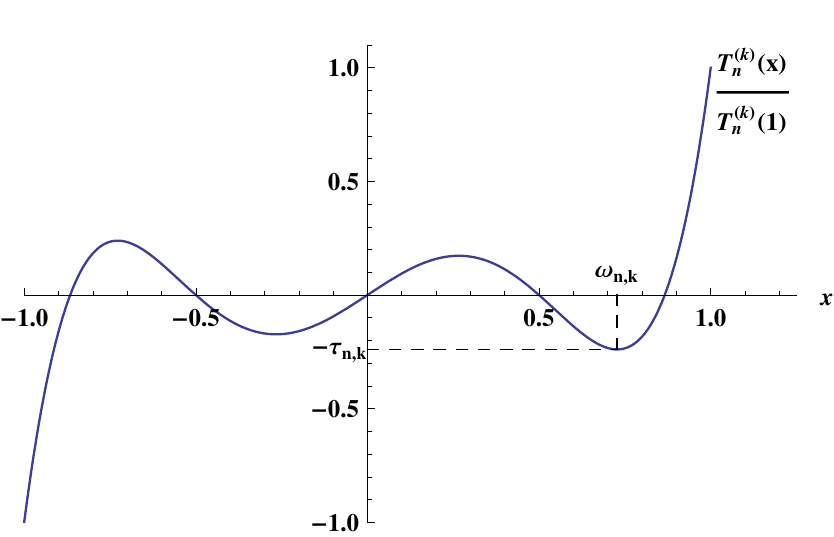}
\caption{{\small The last relative extremum $\,\tau_{n,k}\,$ (here,
$k=1$, $\,n=6$)\,.}}
    \label{figure:1}
\end{figure}

Let us mention previous results.
For the first derivative ($k=1$), Erd\H{o}s--Szeg\H{o} \cite{es} showed
that
$\tau_{3,1} = \frac{1}{3}$, $\tau_{4,1} = \frac{1}{3}(\frac{2}{3})^{1/2}$
and proved that
\be \label{e0}
    \tau_{n,1} \le \frac{1}{4}\,, \qquad n \ge 5\,.
\ee
For arbitrary $k \ge 1$, Eriksson \cite{e} and Nikolov
\cite{n1} independently showed that
\be \label{e1}
    \tau_{n,k} \le \frac{1}{2k+1}\,, \qquad n \ge k+2\,,
\ee
with a better estimate when $n$ is large relative to $k$,
\be \label{e2}
   \tau_{n,k} \le \frac{1}{2k+1}\,\frac{8}{2k+7}\,, \qquad
   n \gtrsim k^{3/2}.
\ee (The exact condition in \cite{e}, \cite{n1} was $\omega_{n,k}
\ge 1 - \frac{8}{2k+7}$ which implies the above inequality between
$n$ and $k$ via the upper estimate $\omega_{n,k} <  1 -
\frac{k^2}{n^2}$.)

\medskip
In this paper, motivated by the applications mentioned above, we
refine and extend inequalities \rf[e0]-\rf[e2] in several
directions.

\medskip
1) Our first observation is a monotone behaviour of the value
$\tau_{n,k}$ with respect to $n$.

\begin{theorem}\label{t1.1}
For a fixed $k \in \N$, the values $\tau_{n,k}$ decrease monotonically in $n$,
i.e.,
\be \lb{tau<<}
   \tau_{n+1,k} < \tau_{n,k} < \cdots < \tau_{k+3,k} < \tau_{k+2,k}\,.
\ee
In particular, for any fixed $k\in \mathbb{N}\,$ and any $m \ge 2$,
we have
\be \lb{tau<}
   \tau_{n,k} \le \tau_{k+m,k}, \qquad n \ge k + m\,.
\ee
\end{theorem}

\begin{figure}[h]
\centering
\includegraphics[scale=0.9]{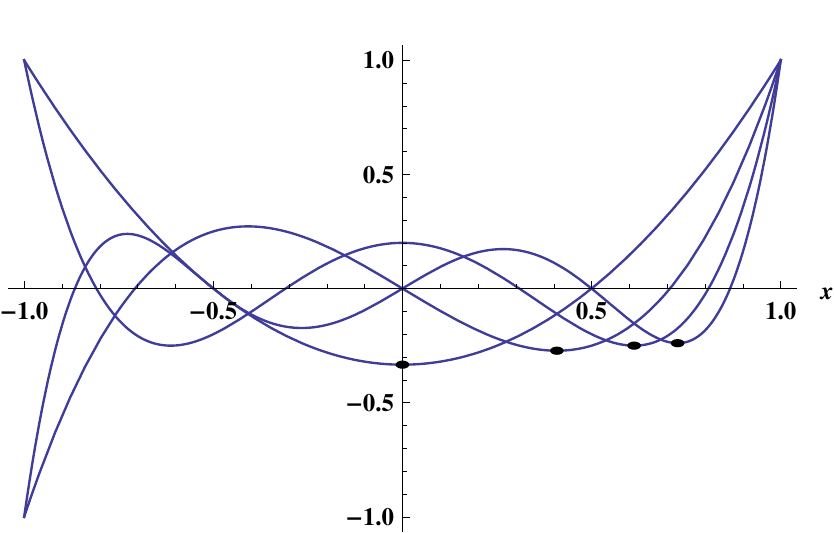}
\caption{{\small The last relative extrema $\,\tau_{n,1}\,$,
$\,3\leq n\leq 6$\,.}}
    \label{figure:2}
\end{figure}

In fact, such a monotone decrease of the relative values of the
local extrema takes place for the ultraspherical polynomials
$P_n^{(\la)}$ with any parameter $\la
> 0$. This remarkable result is due to Sz\'{a}sz \cite{os} and,
for reader's convenience and to keep the paper self-contained, we
state it as Theorem~\ref{t2.1} and give a short proof.

\medskip
2) Our next result is several sharp estimates for $\tau_{n,k}$
which follow from \rf[tau<].
Namely, since $T_{k+m}^{(k)}$ is a symmetric polynomial of degree $m$, for
small $m = 2..6$ we compute the value of its largest extremum, hence
$\tau_{k+m,k}$, explicitly and then use \rf[tau<].

\begin{theorem}\label{tm=4}
We have
\be \lb{m=4}
    \tau_{n,k} \le \left\{
\begin{array}{rcll}
   \tau_{k+2,k}
&=& \frac{1}{2k+1}\,, & n \ge k+2\,, \vspace*{2mm}\\
   \tau_{k+3,k}
&=& \frac{1}{2k+1}\, \big(\frac{2}{k+2}\big)^{1/2}\,,
    & n \ge k+3\,, \vspace*{2mm}\\
   \tau_{k+4,k}
&=& \frac{1}{2k+1}\,\frac{3}{k+3}\,, & n \ge k+4\,.
\end{array} \right.
\ee
\end{theorem}

These estimates contain earlier results \rf[e0]-\rf[e1] as
particular cases, and the last inequality in \rf[m=4] improves
\rf[e2] by the factor of $\frac{3}{4}$ and removes the unnecessary
restrictions on $n$ and $k$.

The next pair of estimates strengthens \rf[m=4]. It also shows
that, although the nice pattern for $\tau_{k+m,k}$ in \rf[m=4] is no
longer true for $m \ge 5$, an approximate behaviour $\tau_{k+m,k}
\approx (\frac{m}{k+m})^{m/2}$ is very much suggestive.

\begin{theorem} \lb{tm=6}
We have
\be\label{m=6}
    \tau_{n,k} \le \left\{
    \begin{array}{rcll}
    \tau_{k+5,k}
&=& \frac{1}{2k+1}\, \big(\frac{4}{k+4}\big)^{3/2} \alpha_k\,,
    & n \ge k+5\,,
    \vspace*{2mm}\\
    \tau_{k+6,k}
&=& \frac{1}{2k+1}\,\big(\frac{5}{k+5}\big)^2 \beta_k\,,
   & n \ge k+6\,,
\end{array} \right.
\ee
where the values $\alpha_k,\,\beta_k$ increase monotonically
to the following limits,
$$
    \alpha_k < \alpha_* = \frac{\sqrt{3(3+\sqrt{6})}}{4}
    = 1.0108.., \qquad
    \beta_k < \beta_* = \frac{2+\sqrt{10}}{5} = 1.0325..
$$
\end{theorem}

Let us note that, because of monotonicity of $\,\tau_{n,k}$,
for any fixed moderate $k$ and a moderate $n_0$, one can compute numerically
the value $\,\tau_{n_0,k}$, thus getting for particular $k$ the estimate
$$
   \tau_{n,k} \le \tau_{n_0,k}\,,\quad n\geq n_0\,.
$$
which would be better than those in \rf[m=4] and \rf[m=6].

\medskip
3) An approximate behaviour
$\tau_{k+m,k} \approx (\frac{m}{k+m})^{m/2}$ in \rf[m=4]-\rf[m=6]
suggests
that when $\,m\,$ is fixed and $k$ grows,
then $\tau_{n,n-m} = \tau_{k+m,k}\,$ is of a
polynomial decay in $n$, i.e.,
$$
    \tau_{n,n-m} = \OO(n^{-m/2}) \qquad (n \to \infty)\,,
$$
while when $k$ is fixed and $n$ grows, we have and exponential estimate
in $k$,
$$
   \tau_{n,k} = \OO(e^{-\gamma k}) \qquad (n \to \infty)\,.
$$
We prove that such a behaviour is indeed the case by
establishing first the upper bounds for $\tau_{n,k}$ which are uniform
in $n$ and $k$, and then considering different relations between
$n$ and $k$.

\begin{theorem}\label{t1.2}
For every $\,n,k\in \N\,$ with $\,n \ge k+2$, we have
\ba
      \tau_{n,k}^2
&\le& \frac{1}{2}
      \Big(1 + \frac{k}{n}\Big)
      \Big(\frac{n}{k}\Big)^{2k}
      {n+k \choose n-k}^{-1} \\
&\le& c_1^2\, k^{\frac{1}{2}}\,\Big(1 - \frac{k^2}{n^2}\Big)^{\frac{1}{2}}
       \frac{ (2n)^{2k} (n-k)^{n-k} } {(n+k)^{n+k}},
       \qquad c_1^2 = \frac{e^2}{2\sqrt{\pi}}\,.
\ea
\end{theorem}

As a consequence of Theorem~\ref{t1.2} we obtain the following statement.

\begin{theorem}\label{t1.3}
We have the following estimates:

\smallskip
(i) if $k \in \N$ is fixed and $n$ grows, then
\be
    \tau_{n,k}
\le c_1 \Big(\frac{2}{e}\Big)^k
        \frac{k^{1/4}}{(1-\frac{k^2}{n^2})^{k/2}}\,,
\ee
in particular
\be
\lb{taunk}
     \tau_{n,k} < c_2 \Big(\frac{2}{e}\Big)^k k^{1/4}, \qquad
     n \ge k^{3/2}\,;
\ee

(ii) if $n-k = m \in \N$ is fixed and $n$ grows, then
\be \lb{k+m}
    \tau_{n,n-m}
\le c_3\, m^{1/4}
    \Big(\frac{me}{2}\Big)^{m/2} n^{-m/2}\,;
\ee

(iii) if $k = \lfloor\la n \rfloor \in \N$, where $\la \in (0,1)$,
and $n$ grows, then we have an exponential decay
$$
     \tau_{n,\la n}
\le c_4\, n^{1/4} \rho_\la^{n/2}, \qquad \rho_\la < 1\,,
$$
in particular
$$
   \tau_{n,n/2} < c_1 n^{1/4} \Big(\frac{4}{\sqrt{27}}\Big)^{n/2}.
$$
\end{theorem}

We can reformulate Theorem \ref{t1.3} in the form which shows,
for a fixed $k$ and growin $n$,
the rate of decrease of the values $\tau_{n,k}$ in \rf[tau<<].

\begin{corollary}
We have
\be
    \tau_{n,k}
\lesssim \left\{ \begin{array}{cl}
    k^{-m/2}, & n \ge k+m\,, \\
    (\frac{4}{\sqrt{27}})^k, & n \ge 2k\,, \\
    (\frac{2}{e})^k, & n \ge k^{3/2}\,.
    \end{array} \right.
\ee
\end{corollary}

\begin{remark} \rm
Exponential estimate \rf[taunk] becomes superior to the polynomial
esimates \rf[m=6] only when $k \ge 10$.
\end{remark}

4) Finally, we establish the asymptotics of the values
of $\lim\limits_{n\to \infty}\tau_{n,k}\,$ and
$\lim\limits_{n\to \infty}\tau_{n,n-m}$
which shows that
the upper bounds in \rf[taunk] and \rf[k+m] are
asymptotically correct with respect to the exponential terms therein.

\begin{theorem}\label{t1.4}
We have
\ba
    \tau_k^* := \lim_{n\to\infty} \tau_{n,k}
&=& C_0 \Big(\frac{2}{e}\Big)^k e^{-a_0 k^{1/3}}
   k^{-1/6}\big(1+\OO(k^{-1/3})\big), \lb{tau*}\\
  \tau_m^{**} := \lim_{n\to\infty} n^{m/2}\tau_{n,n-m}
&=& C_1
  \Big(\frac{e\,m}{2}\Big)^{m/2} e^{-a_1 m^{1/3}}
  m^{-1/6}\big(1+\OO(m^{-1/3})\big)\,, \lb{tau**}
\ea
where the pairs of constants
\baa
   a_0 = 1.8557...\,, \quad C_0 = 1.1966..
   \quad\mbox{and}\quad
   a_1 = 2.3381...\,, \quad C_1 = 1.0660...
\eaa
can be explicitly represented in terms of the Airy function.
\end{theorem}

\medskip
The rest of the paper is organised as follows.
In Section~2, we present a proof of Theorem~\ref{t1.1}, which is
deduced from a more general statement, Theorem~\ref{t2.1}, about
monotonicity of the relative extrema of ultraspherical polynomials.
In Section \ref{smallm},  we compute directly $\tau_{k+m,k}$ for
small $m$, thus proving Theorems~\ref{tm=4}-\ref{tm=6}.
In Section~4, we adopt the majorant, originally introduced by
Shaeffer and Duffin \cite{sd} for their alternative proof of the
Markov inequality, and prove then Theorem~\ref{t1.2}.
Theorem~\ref{t1.3} is proved in Section~5.
The proof of Theorem~\ref{t1.4}, which relies on some known
asymptotic behaviour of orthogonal polynomials, is given in
Section~\ref{as}.


\section{Monotonicity of the sequence
\boldmath{$\{\tau_{n,k}\}_{n\geq k+2}$}}


Here, we prove that
$$
    \mu_{i,n}^{(\la)} = \frac{P_n^{(\la)}(y_{i,n}^{(\la)})}{P_n^{(\la)}(1)}\,,
$$
the relative values of the ordered local extrema of the
ultraspherical polynomials $P_n^{(\la)}$ with parameter $\la$ decay
monotonically with respect to $n$ for any $\la > 0$. This includes
Theorem \ref{t1.1} as a particular case since $T_n^{(k)}$, the
$k$-th derivative of the Chebyshev polynomials of degree $n$,
coincide up to a factor with $P_{n-k}^{(\la)}$, where $\la = k$.

We start with recalling some known facts about the ultraspherical
polynomials (for more details, see \cite[Chapther 4.7]{GS}).

For $\,\la> -\frac{1}{2}$, $\{P_n^{(\la)}\}_{n\in \N_0}$ stands for
the sequence of ultraspherical polynomials, which are orthogonal on
$\,[-1,1]\,$ with respect to the weight function
$w_{\la}(x)=(1-x^2)^{\la-1/2}$, with the standard normalization
$$
    P_n^{(\la)}(1) = {n+2\la-1\choose n}\,,\qquad \la\ne 0.
$$
The Chebyshev polynomials of the
first and the second kind and the Legendre polynomials are
particular cases of ultraspherical polynomials, they correspond
up to a factor to
the values $\,\la=0,\,1\,$ and $\frac{1}{2}$, respectively.
Moreover, due to the properties
\baa
     T_n^{\prime}(x)
 &=& n\,P_{n-1}^{(1)}(x)\,,\vspace*{2mm}\\
     \frac{d}{dx}\, P_{n}^{(\la)}(x)
 &=&  2\la\,P_{n-1}^{(\la+1)}(x)\,, \qquad \la\ne 0\,,
\eaa
the derivatives of the Chebyshev polynomials are ultraspherical
polynomials, too,
\be  \label{e2.1}
    T_n^{(k)}(x)
  = c_{n,\la} P_{n-k}^{(\la)}(x)\,,\qquad \la = k, \qquad
    k=1,\ldots,n\,.
\ee

We will work with the re-normalised ultraspherical polynomials
\be \label{e2.2}
    p_n^{(\la)}(x)
 := P_n^{(\la)}(x)/P_n^{(\la)}(1),
\ee
so that $p_n^{(\la)}(1)=1$. It is clear that the absolute
values of the local extrema of $p_n^{(\la)}$ are equal
to the relative values of the local extrema of $P_n^{(\la)}$
compared to $P_n^{(\la)}(1)$.

Theorem~\ref{t1.1} is a consequence of the following general statement.

\begin{theorem}\label{t2.1}
Let $y_{1,n}^{(\la)} > y_{2,n}^{(\la)} > \cdots > y_{n-1,n}^{(\la)}$
be the zeros of the ultraspherical polynomial $p_{n-1}^{(\la+1)}$,
i.e., the abscissae of the local extrema of $p_n^{(\la)}$, in the
reverse order. Set $y_{n,n}^{(\la)} := -1$, and denote
$$
    \mu_{i,n}^{(\la)}
 := |p_n^{(\la)}(y_{i,n}^{(\la)})|, \qquad
    i=1,\ldots,n\,.
$$

1) If $\lambda > 0$, then
\be \label{e2.5}
   \mu_{i,n+1}^{(\la)} < \mu_{i,n}^{(\la)} \quad \mbox{ for }\quad
i=1,2,\ldots,n\,.
\ee

2) If $-\frac{1}{2} < \la <0$, then inequalities \rf[e2.5] hold with the
opposite sign.

3) (If $\la=0$, then $p_n^{(\la)} = T_n$ and
we have equalities in \rf[e2.5]
as all local extrema of $T_n$ and $T_{n+1}$ are of the absolute value $1$.)
\end{theorem}

\proof We omit index $\la$, so set $p_n := p_n^{(\la)}$, and we will
use the next two identities which readily follow from
\cite[eqn.(4.7.28)]{GS}: \baa
    p_n(x)
&=& - \frac{1}{n+2\la}\,x\,p_n^{\prime}(x)
    + \frac{1}{n+1}\,p_{n+1}^{\prime}(x)\,, \\
    p_{n+1}(x)
&=& - \frac{1}{n+2\la}\,p_n^{\prime}(x)
    + \frac{1}{n+1}\,x\,p_{n+1}^{\prime}(x)\,.
\eaa
From those we deduce that
\be \lb{p}
     p_n(x)^2 - p_{n+1}(x)^2
  =  (1-x^2)\Big[\frac{1}{(n+1)^2}\,p_{n+1}^{\prime}(x)^2
        -\frac{1}{(n+2\lambda)^2}\,p_n^{\prime}(x)^2\Big]\,,
\ee
and we rearrange his equality as follows,
\be \label{e2.6}
   f(x)
 := p_n(x)^2+\frac{1-x^2}{(n+2\lambda)^2}\,p_n^{\prime}(x)^2
\stackrel{\rf[p]}{=}
   p_{n+1}(x)^2+\frac{1-x^2}{(n+1)^2}\,p_{n+1}^{\prime}(x)^2\,.
\ee
Clearly, $\,f\,$ is a polynomial of degree $2n$ which interpolates
both $\,p_n^2\,$ and $\,p_{n+1}^2\,$ at the points of their local
maxima in $\,[-1,1]$. Moreover, $\,f^{\prime}\,$ vanishes at the
zeros of both $\,p_n^{\prime}\,$ and $\,p_{n+1}^{\prime}\,$,
therefore, with some constant $\,c_n^{(\la)}$,
\be \label{e2.7}
    f^{\prime}(x)
  = c_n^{(\la)}\,p_n^{\prime}(x)p_{n+1}^{\prime}(x)\,.
\ee

Next, we determine the sign of $\,c_n^{(\la)}$. Let $a_n, a_{n+1}$ be
the leading coefficients of $p_n$ and $p_{n+1}$, respectively,
and note that, since $p_n(1) = p_{n+1}(1) = 1$, we have $a_n, a_{n+1} > 0$.
Then, equating the leading coefficients
of $f'$ in representations \rf[e2.6] and \rf[e2.7], respectively,
we obtain
$$
   2n a_n^2 \Big(1-\frac{n^2}{(n+2\la)^2}\Big)
 = c_n^{(\la)} n(n+1)\, a_n a_{n+1}\,,
$$
whence
$$
   c_n^{(\la)}
 = \frac{1}{n+1}\frac{a_n}{a_{n+1}}\frac{4\la (2n+2\la)}{(n+2\la)^2}
 \ra \sign c_n^{(\la)} = \sign \la
$$
Thus, \rf[e2.7] becomes
\be \label{e2.8}
   f^{\prime}(x)
 = c\, \la\, p_n^{\prime}(x)\,p_{n+1}^{\prime}(x)\,, \qquad
   c = c_{n,\la} > 0\,.
\ee

Now we can prove Theorem~\ref{t2.1}.
Let $\la > 0$. Then from
\rf[e2.8] and the interlacing of zeros of $p_n^{\prime}$ and
$p_{n+1}^{\prime}$ we conclude that
$$
    f^{\prime}(x) < 0\,,\quad
    x\in (y_{i,n}^{(\la)},y_{i,n+1}^{(\la)})\,,\quad i=1,\ldots,n\,,
$$
i.e., $\,f\,$ is monotonically decreasing on each interval
$\,(y_{i,n}^{(\la)},y_{i,n+1}^{(\la)})$. From \rf[e2.6], we have
\baa
     f(y_{i,n}^{(\la)})
 &=& p_{n}(y_{i,n}^{(\la)})^2 \;=\; |\mu_{i,n}^{(\la)}|^2\,, \\
     f(y_{i,n+1}^{(\la)})
 &=& p_{n+1}(y_{i,n+1}^{(\la)})^2 \;=\; |\mu_{i,n+1}^{(\la)}|^2\,,
\eaa
therefore
$$
    |\mu_{i,n+1}^{(\la)}| < |\mu_{i,n}^{(\la)}|, \qquad i =1,\ldots,n\,.
$$
Clearly, if $\la < 0$, then the sign is reversed.
\qed

\medskip
{\bf Proof of Theorem~\ref{t1.1}.}
By \rf[e2.1] and \rf[e2.2], we have
$$
   \frac{T_n^{(k)}(x)}{T_n^{(k)}(1)}
 = \frac{P_{n-k}^{(\la)}(x)}{P_{n-k}^{(\la)}(1)}
 = p_{n-k}^{(\la)}(x)\,, \qquad \la = k \ge 1.
$$
Hence, $\,\tau_{n,k}=\mu_{1,n-k}^{(\la)}\,$,
and then Theorem~\ref{t2.1} yields
$$
   \tau_{n,k} = \mu_{1,n-k}^{(\la)}
 > \mu_{1,n+1-k}^{(k)} = \tau_{n+1,k}\,.
$$
Theorem~\ref{t1.1} is proved. \qed
\medskip


\section{Proof of Theorems \ref{tm=4}-\ref{tm=6} } \lb{smallm}


By Theorem~\ref{t1.1}, for any fixed $\,m$, the value
$\,\tau_{k+m,k}\,$ gives an upper bound for all $\,\tau_{n,k}\,$,
namely
$$
    \tau_{n,k} \le \tau_{k+m,k}\,, \qquad n \ge k+m,
$$
so here we determine the latter values directly for $m = 2..6$.

\smallskip
We will need the expansion formula for the $n$-th Chebyshev
polynomial,
\ba
    T_n(x)
&=& \frac{n}{2} \sum_{i=0}^{\lfloor n/2 \rfloor}
    (-1)^i \frac{(n-i-1)!}{i!(n-2i)!} (2x)^{n-2i} \nonumber \\
&=& 2^{n-1}\, x^n \;-\; 2^{n-3}\, n\, x^{n-2} \\
& &   \; +\;  2^{n-6}\,n(n-3)\, x^{n-4}
  \;-\; \frac{1}{3}\, 2^{n-8}\, n(n-4)(n-5)\, x^{n-6} + \; \cdots
    \lb{e2.9}
\ea
From this we compute expression for $T_n^{(n-6)}$ in \rf[n-6],
and then differentiate it to find
all further derivatives $T_n^{(n-m)}$ for $m=5..2$.

We will denote the point of the rightmost local extrema of $T_n^{(k)}$
by $x_*$, i.e., $x_* := \omega_{n,k}$. Since $T_n^{(k+1)}(x_*) = 0$ then,
for the value of $T_n^{(k)}(x_*)$
we will also use simplifications arising from the formula
$$
   T_n^{(k)}(x_*) =  T_n^{(k)}(x_*) - c_{n,k} x_* T_n^{(k+1)}(x_*)
$$
where we choose the constant $c_{n,k}$ to cancel high degree monomials.

\medskip
{\bf\boldmath 1) The case $k = n-2$ (or equivalently $n = k+2$).}
We have
\be \lb{n-2}
    T_n^{(n-2)}(x)
  = c^{-1}\,\Big[2(n-1)x^2 -1\Big]\,,
\ee
whence $x_* = 0$ and
$$
    \tau_{n,n-2}
  = \frac{|T_n^{(n-2)}(x_*)|}{ T_n^{(n-2)}(1)}
  = \frac{1}{2n - 3}
\ra \tau_{k+2,k} = \frac{1}{2k+1}\,.
$$

\medskip
{\bf\boldmath  2) The case $k = n-3$ (or equivalently $n = k+3$).}
We obtain
\be \lb{n-3}
   T_n^{(n-3)}(x)
= c^{-1}\, \Big[2(n-1)x^3 -3x\Big]\,,
\ee
hence  $c\, T_n^{(n-3)}(1) = 2n-5$. From \rf[n-2], we find
$x_*^2 = \frac{1}{2(n-1)}$ and
$$
    c\,T_n^{(n-3)}(x_*) = - 2x_* = -\frac{2}{\sqrt{2(n-1)}}\,.
$$
Respectively,
$$
   \tau_{n,n-3}
=  \frac{|T_n^{(n-3)}(x_*)|}{ T_n^{(n-3)}(1)}
=  \frac{1}{2n-5}\sqrt{\frac{2}{n-1}}
\ra \tau_{k+3,k} = \frac{1}{2k+1}\,\sqrt{\frac{2}{k+2}}\,.
$$

\medskip
{\bf\boldmath 3) The case $k = n-4$ (or equivalently $n = k+4$).}
We have
\be \lb{n-4}
   T_n^{(n-4)}(x)
 = c^{-1}\, \Big[4 (n-1)(n-2) x^4
        - 12 (n-2)x^2 + 3 \Big]\,.
\ee
hence
$$
   c\,T_n^{(n-4)}(1) = 4 n^2 - 24 n  + 35 = (2n-5)(2n-7)\,.
$$
From \rf[n-3], we find $x_*^2 = \frac{3}{2(n-1)}$
and
$$
    c\,T_n^{(n-4)}(x_*) = - 6(n-2)x_*^2 + 3
  = \frac{3}{n-1} \Big[-3(n-2) + (n-1)\Big]
  = - \frac{3(2n-5)}{n-1}\,.
$$
Respectively,
$$
     \tau_{n,n-4}
  =  \frac{|T_n^{(n-4)}(x_*)|}{ T_n^{(n-4)}(1)}
  =  \frac{1}{2n-7}\, \frac{3}{n-1}
 \ra \tau_{k+4,k} = \frac{1}{2k+1}\, \frac{3}{k+3}\,.
$$

The cases 1)-3) prove estimates \rf[m=4], hence Theorem \ref{tm=4}.
\qed

\medskip
{\bf\boldmath 4) The case $k = n-5$ (or equivalently $n = k+5$).}
We have
\be \lb{n-5}
   T_n^{(n-5)}(x)
 = c^{-1} \Big[ 4(n-1)(n-2) x^5 -  20 (n-2)x^3 + 15x \Big],
\ee
hence
$$
   c\,T_n^{(n-5)}(1) =  4n^2 - 32n + 63 =  (2n-9)(2n-7)\,.
$$
From \rf[n-4], we find
\baa
   x_*^2
&=& \frac{3(n-2) + \sqrt{9(n-2)^2 - 3(n-1)(n-2)}}{2(n-1)(n-2)} \\
&=& \frac{1}{n-1} \frac{3 + \sqrt{6-t}}{2}, \qquad
   t := t_k = \frac{3}{n-2} = \frac{3}{k+3}\,.
\eaa
and
$$
   c\,T_n^{(n-5)}(x_*) = -4x_* \Big[2(n-2)x_*^2-3 \Big].
$$
After simplifications we obtain
$$
   \tau_{n,n-5}
 = \frac{|T_n^{(n-5)}(x_*)|}{T_n^{(n-5)}(1)}
 = \frac{1}{2n-9}\frac{8}{(n-1)^{3/2}}\, \alpha_k
\ra \tau_{k+5,k}
   = \frac{1}{2k+1}\frac{4^{3/2}}{(k+4)^{3/2}}\, \alpha_k\,,
$$
where
$$
   \alpha_k
:= \frac{1}{2\sqrt{2}}
   \frac{\sqrt{3 + \sqrt{6-t}}}{2-t}\Big(\sqrt{6-t} -t\Big)
 = \frac{1}{2\sqrt{2}}
   \frac{(y+3)^{3/2}}{y+2} =: f(y), \qquad y := \sqrt{6-t}\,.
$$
The function $f$ is increasing for $y > 0$, hence
$$
   t_k > t_{k+1} \ra y_k < y_{k+1} \ra \alpha_k < \alpha_{k+1} < \alpha_*
$$
where
$$
   \alpha_*
 = \lim_{t \to 0} \alpha_k
 = \frac{1}{2\sqrt{2}} \frac{ \sqrt{3+\sqrt{6}}}{2} \sqrt{6}
 = \frac{\sqrt{ 3(3+\sqrt{6})}}{4}\,.
$$

\medskip
{\bf\boldmath 5) The case $k = n-6$ (or equivalently $n = k+6$).}
From \rf[e2.9], we have
\be \lb{n-6}
   T_n^{(n-6)}(x)
 = c^{-1} \Big[ 8(n-1)(n-2)(n-3) x^6 -  60 (n-2)(n-3)x^4
    + 90(n-3) x^2 - 15 \Big],
\ee
hence
$$
   c\,T_n^{(n-6)}(1)
 =  8n^3 - 108n^2 + 478n - 693 \\
=  (2n-11)(2n-9)(2n-7)\,.
$$
From \rf[n-5], we find
\baa
   x_*^2
&=& \frac{5(n-2) + \sqrt{25(n-2)^2 - 15(n-1)(n-2)}}{2(n-1)(n-2)} \\
&=& \frac{1}{n-1} \frac{5 + \sqrt{10-3t}}{2}, \qquad
   t := t_k = \frac{5}{n-2} = \frac{5}{k+4}
\eaa
and
$$
   c\,T_n^{(n-6)}(x_*) = -4(2n-7)\, x_*^2 \Big[2(n-2)x_*^2-5 \Big]\,.
$$
After simplifications, we obtain
$$
   \tau_{n,n-6}
 = \frac{|T_n^{(n-6)}(x_*)|}{T_n^{(n-6)}(1)}
 = \frac{1}{2n-11}\frac{5^2}{(n-1)^2}\, \beta_k
\ra \tau_{k+6,k} = \frac{1}{2k+1}\frac{5^2}{(k+5)^2}\, \beta_k\,,
$$
where
$$
   \beta_k
:= \frac{2}{5^2}
   \frac{5 + \sqrt{10-3t}}{2-t}\Big(\sqrt{10-3t} -t\Big)
 = \frac{2}{5^2}
   \frac{(y+5)^2}{y+2} =: g(y), \qquad y := \sqrt{10-3t}\,.
$$
The function $g$ is increasing for $y > 1$, hence
$$
    t_k > t_{k+1} \ra y_k < y_{k+1}
\ra \beta_k < \beta_{k+1} < \beta_*\,,
$$
where
$$
   \beta_*
 = \lim_{t \to 0} \beta_k
 = \frac{2}{5^2} \frac{5+\sqrt{10}}{2} \sqrt{10}
 = \frac{2+\sqrt{10}}{5}\,.
$$

The cases 4)-5) prove estimates \rf[m=6], hence Theorem \ref{tm=6}.
\qed


\section{Estimates based on the Duffin--Shaeffer majorant}


In this section, we prove Theorem~\ref{t1.2}.
Our proof is based on the upper bound $\tau_{n,k} < \delta_{n,k}$
which uses the so-called Duffin--Schaeffer majorant.

\begin{definition} \rm
With $T_n$ the Chebyshev polynomial of degree $n$, and $S_n(x) :=
\frac{1}{n} \sqrt{1-x^2}\,T_n'(x)$, we define the
{\it Duffin--Schaeffer majorant} $D_{n,k}(\cdot)$ as
\begin{equation} \label{Ddef}
    D_{n,k}(x) :=  \{[T_n^{(k)}(x)]^2 + [S_n^{(k)}(x)]^2\}^{1/2},
   \quad x \in (-1,1).
\end{equation}
\end{definition}

This majorant was introduced by Shaeffer--Duffin \cite{sd} who
proved that, if $\,p\,$ is a polynomial of degree not exceeding $n$,
then
\be \label{sd1}
    \|p\| \le 1
\ra |p^{(k)}(x)| \le D_{k,n}(x)\,, \quad x \in (-1,1).
\ee
which may be viewed as a generalization of the pointwise
Bernstein inequality $|p'(x)| \le \frac{n}{\sqrt{1-x^2}} \|p\|$
to higher derivatives.

\begin{lemma}
The majorant $\,D_{n,k}\,$ has the following properties.
\begin{enumerate}
\item
We have
\begin{equation} \label{TD}
|T_n^{(k)}(x)| \le D_{n,k}(x) \quad\mbox{ for all } x \in (-1,1)\,.
\end{equation}

\item
$\,D_{n,k}(x) = |T_n^{(k)}(x)|$ at zeros of $S_n^{(k)}$, in
particular,
\begin{equation} \label{D0}
   D_{n,k}(0) = |T_n^{(k)}(0)| \quad \mbox{ if $n-k$ is even.}
\end{equation}

\item
The majorant $\,D_{n,k}(\cdot)\,$ is a strictly increasing function
on $\,[0,1]$.

\item
We have the explicit formulae
$\,\displaystyle{\frac{1}{n^2}\,\big[D_{n,1}(x)\big]^2=\frac{1}{1-x^2}}\,$
and
\begin{equation} \label{Df}
   \frac{1}{n^2}\,\big[D_{n,k}(x)\big]^2
= \sum_{m=0}^{k-1} \frac{b_{m,n}}{(1-x^2)^{k+m}}\,,\qquad k\geq 2\,,
\end{equation}
where
\ba \label{b_m}
     b_{m,n}
 &=& c_{m,k} (n^2-(m+1)^2)\cdots(n^2-(k-1)^2)\,, \\
     c_{m,k}
 &:=& \left\{\begin{array}{ll} 1, & m=0\,, \\
      {k-1+m \choose 2m}(2m-1)!!^2\,, & m\geq 1\,.
\end{array}\right.
\ea
\end{enumerate}
\end{lemma}

\proof
Claim 1 and the first half of Claim 2 follow
directly from Definition \ref{Ddef}. Equality \rf[D0] is due to the
fact that $T_n$ and $S_n$ are of different parity, so if $n-k$
is even, then $T_n^{(k)}$ is an even function and $S_n^{(k)}$ is an
odd one, hence $S_n^{(k)}(0)=0$. The third property was proved by
Schaeffer--Duffin \cite{sd}, and it also follows easily from the
formulas \rf[Df]-\rf[b_m] which were established by Shadrin
\cite{s1}. \qed

\smallskip
Here are few particular expressions for $D_{n,k}(\cdot)$.
\baa
    \frac{1}{n^2}[D_{n,1}(x)]^2
& = & \displaystyle{
       \frac{1}{1-x^2} }\,, \\
    \frac{1}{n^2}[D_{n,2}(x)]^2
& = & \displaystyle{
      \frac{(n^2-1)}{(1-x^2)^2} + \frac{1}{(1-x^2)^3}}\,,\\
    \frac{1}{n^2} [D_{n,3}(x)]^2
& =& \displaystyle{
      \frac{(n^2-1)(n^2-4)}{(1-x^2)^3} + \frac{3(n^2-4)}{(1-x^2)^4}
       + \frac{9}{(1-x^2)^5}}\,,\\
    \frac{1}{n^2} [D_{n,4}(x)]^2
& = & \displaystyle{
      \frac{(n^2-1)(n^2-4)(n^2-9)}{(1-x^2)^4} +
        \frac{6(n^2-4)(n^2-9)}{(1-x^2)^5}} \\
   && \displaystyle{
       \quad +\, \frac{45(n^2-9)}{(1-x^2)^6} + \frac{225}{(1-x^2)^7}}\,.
\eaa

\begin{lemma}\label{l3.2}
Let $\,\omega_k :=\omega_{n,k}\,$ be the rightmost zero of
$\,T_n^{(k+1)}$. Then
\be \label{xk}
    \omega_k < x_k, \qquad \mbox{where}\quad
    x_k^2 := 1 - \frac{k^2}{n^2}\,.
\ee
\end{lemma}

\proof
The claim can be deduced from numerous upper bounds for
the extreme zeros of ultraspherical polynomials. For instance, in
\cite{n1} Nikolov proved that $\omega_k^2 \le \frac{n^2 -
(k+2)^2}{n^2 + \alpha_{n,k}}$, with some $\alpha_{n,k} > 0$, hence
\be \lb{w<}
    \omega_k^2
\le \frac{n^2 - (k+2)^2}{n^2}
\le \frac{n^2 - k^2}{n^2} = x_k^2\,.
\ee
\qed

From \rf[TD], monotonicity of $D_{n,k}(\cdot)$ and inequality \rf[xk],
it follows immediately that
$$
     |T_n^{(k)}(\omega_k)| \le D_{n,k}(\omega_k) < D_{n,k}(x_k),
$$
hence the following statement.

\begin{proposition}\label{p3.1}
We have
$$
    \tau_{n,k} < \delta_{n,k}, \qquad
    \delta_{n,k} := \frac{D_{n,k}(x_k)}{T_n^{(k)}(1)}\,.
$$
\end{proposition}

We proceed with estimates of $\,\delta_{n,k}\,$, using the
explicit expression \rf[Df] for $\,D_{n,k}(\cdot)$.

\begin{lemma}\label{l3.3}
We have
\be \label{AB}
    \tau_{n,k}^2 < \delta_{n,k}^2 = A_{n,k} B_{n,k}\,,
\ee
where
\ba
     A_{n,k}
 &=& \frac{(2k-1)!!^2}{k^{2k}}
     \sum_{m=0}^{k-1} \frac{c_{m,k}}{k^{2m}}
     \frac{n^{2m}}{(n^2-1^2)\cdots(n^2-m^2)}\,,  \label{A} \\
     B_{n,k}
 &=& \frac{n^{2k}}{ n^2(n^2-1^2)\cdots(n^2-(k-1)^2)}\,. \label{B}
\ea
\end{lemma}

\proof
From \rf[Df] -- \rf[b_m], we obtain
\begin{eqnarray*}
    [\delta_{n,k}]^2
&=& \frac{[D_{k,n}(x_k)]^2}{[T_n^{(k)}(1)]^2} \\
&=& \frac{1}{[T_n^{(k)}(1)]^2}\, n^2 \sum_{m=0}^{k-1}
    \frac{c_{m,k}}{(1-x_k^2)^{k+m}}
       (n^2-(m+1)^2)\cdots(n^2-(k-1)^2) \\
&=& \frac{n^2(n^2\!-\!1^2)\cdots(n^2\!-\!(k\!-\!1)^2)}
         {[T_n^{(k)}(1)]^2(1-x_k^2)^k}
     \sum_{m=0}^{k-1}  \frac{c_{m,k}}{(1-x_k^2)^m\,}
         \frac{1}{(n^2\!-\!1^2)\cdots(n^2\!-\!m^2)}
\end{eqnarray*}
and substitution
$$
   \frac{1}{[T_n^{(k)}(1)]^2}
 = \frac{(2k-1)!!^2}{[n^2(n^2-1^2)\cdots(n^2-(k-1)^2)]^2}\,, \qquad
   1 - x_k^2 = \frac{k^2}{n^2}\,,
$$
gives \rf[AB] -- \rf[B] after a rearrangement. \qed

\begin{remark} \rm
Whereas the value $\tau_{n,k}$ is defined only for $n \ge k+2$,
the values of $A_{n,k}$ and $B_{n,k}$ in \rf[A]-\rf[B]
are well-defined for $n \ge k$. We will use this fact in the next lemma
where the values $A_{k,k}$ and $B_{k,k}$ will be considered.
\end{remark}

\begin{lemma}\label{ABi}
We have
\be \label{AB0}
    \tau_{n,k}^2 < \delta_{n,k}^2 = A_{n,k} B_{n,k}\,,
\ee
where
\be
    A_{n,k}
\le \frac{1}{2}\, \frac{(2k)!}{k^{2k}}\,, \qquad
    B_{n,k}
=  \frac{n+k}{n} \frac{n^{2k} (n-k)!}{(n+k)!}\,. \lb{AB1}
\ee
\end{lemma}

\proof
Expression for $B_{n,k}$ in \rf[AB1] is just a rearrangment of \rf[B].

As to the inequality for $A_{n,k}$ in \rf[AB1],
it is clear from \rf[A] that $\,A_{n,k}\,$
decreases when $n$ grows, therefore
$$
    A_{n,k} \le A_{k,k}\,, \qquad n \ge k\,.
$$
With $n = k$, we have $\,x_k = 1 - \frac{k^2}{n^2} = 0$, and also
$D_{k,k}(0) = T_{k}^{(k)}(0)$ by \rf[D0], therefore
$$
   A_{k,k} B_{k,k}
 = [\delta_{k,k}]^2
 = \frac{[D_{k,k}(0)]^2}{[T_k^{(k)}(1)]^2}
\stackrel{\rf[D0]}{=} \frac{[T_{k}^{(k)}(0)]^2}{[T_{k}^{(k)}(1)]^2}
 = 1,
$$
hence, $A_{k,k} = 1/B_{k,k}$, and from formula \rf[AB1], we find
$$
   A_{k,k} = \frac{1}{B_{k,k}}
 = \frac{k}{2k}\frac{(2k)!}{k^{2k}}\,,
$$
hence the result.
\qed

\begin{remark} \rm
If we consider the first estimate in \rf[w<], namely
$$
    \omega_k \le x_k', \quad\mbox{where}\quad
    x_k'^2 = 1 - \frac{(k+2)^2}{n^2}, \qquad n \ge k+2\,,
$$
then we obtain
$$
    A_{n,k} \le A'_{k+2,k}
 =  \gamma_k^2 \frac{1}{2} \frac{(2k)!}{k^{2k}}, \qquad
    \gamma_k^2 = \frac{(k+2)}{(2k+1)}\Big(\frac{k}{k+2}\Big)^{2k}
$$
i.e., we can improve the estimate \rf[AB1] (and all subsequent estimates)
by the factor of $\gamma_k$ (or $\gamma_k^2$). Note that
$$
   \gamma_k \approx \frac{1}{\sqrt{2}} \frac{1}{e^2}\,.
$$
\end{remark}

Now, we prove Theorem \ref{t1.2} which is the following statement.

\begin{theorem}\label{t1.2'}
For every $\,n,k\in \N\,$ with $\,n \ge k+2$, we have
\ba
      \tau_{n,k}^2
&\le& \frac{1}{2}
      \Big(1+\frac{k}{n}\Big)
      \Big(\frac{n}{k}\Big)^{2k}
      {n+k \choose n-k}^{-1} \\
&\le& c_1^2 k^{\frac{1}{2}}\,\Big(1-\frac{k^2}{n^2}\Big)^{\frac{1}{2}}
       \frac{ (2n)^{2k} (n-k)^{n-k} } {(n+k)^{n+k}},
       \qquad c_1^2 = \frac{e^2}{2\sqrt{\pi}}\,. \lb{p2}
\ea
\end{theorem}

\proof
The first part is just the estimate \rf[AB0],
$$
   \tau_{n,k}^2 < \delta_{n,k}^2
= A_{n,k} B_{n,k}
< \frac{1}{2} \frac{n+k}{n} \frac{n^{2k}}{k^{2k}}
    \frac{ (n-k)! (2k)!}{ (n+k)!}\,.
$$

To prove the second inequality we use the following version of
Stirling's formula
\be \label{e3.11}
     \sqrt{2\pi} \Big(\frac{N}{e}\Big)^N \sqrt{N}
   < N!
   < e \Big(\frac{N}{e}\Big)^N \sqrt{N}\,.
\ee
This gives
\baa
      \frac{1}{2}\,\frac{n+k}{n} \frac{n^{2k}}{k^{2k}}
      \frac{(n-k)! (2k)!}{(n+k)!}
&\le& \frac{1}{2}\,\frac{n+k}{n} \frac{n^{2k}}{k^{2k}}
      \frac{e^2}{\sqrt{2\pi}} \frac{ \sqrt{n-k}\sqrt{2k}}{\sqrt{n+k}}
      \frac{(n-k)^{n-k} (2k)^{2k}} {(n+k)^{n+k}} \\
& = & \frac{e^2}{2\sqrt{\pi}}
       \frac{n+k}{n} \frac{\sqrt{n-k}\sqrt{k}}{\sqrt{n+k}}
      \frac{ (n-k)^{n-k} (2n)^{2k}} {(n+k)^{n+k}} \\
& = & c_1^2
      k^{\frac{1}{2}}\,\Big(1-\frac{k^2}{n^2}\Big)^{\frac{1}{2}}
      \frac{ (n-k)^{n-k} (2n)^{2k}} {(n+k)^{n+k}}\,,
\eaa
and that finishes the proof.
\qed


\section{Proof of Theorem \ref{t1.3}}


We rewrite inequality \rf[p2] in a more convenient form
\be
  \tau_{n,k}
 \le c_1^2 k^{\frac{1}{2}}\,\Big(1-\frac{k^2}{n^2}\Big)^{\frac{1}{2}}
       \Big(\frac{2n}{n+k}\Big)^{n+k}
       \Big(\frac{n-k}{2n}\Big)^{n-k}\,,
       \qquad c_1^2 = \frac{e}{2\sqrt{\pi}}\,. \lb{p3}
\ee
We will prove each part of Theorem \ref{t1.3} as a separate lemma.

\begin{lemma}
If $k \in \N$ is fixed and $n$ grows, then
\be \lb{tauk}
    \tau_{n,k}
\le c_1 \Big(\frac{2}{e}\Big)^k
        \frac{k^{1/4}}{(1-\frac{k^2}{n^2})^{k/4}}, \qquad
\ee
in particular
\be \lb{tauk1}
     \tau_{n,k} < c_2 \Big(\frac{2}{e}\Big)^k k^{1/4}, \qquad
     n \ge k^{3/2}\,.
\ee
\end{lemma}

\proof
We write \rf[p3] in the form
$$
    \tau_{n,k}^2 \le c_1^2 L_1 L_2\,,
$$
where
\be \lb{es_k}
   L_1 := k^{\frac{1}{2}}\,\Big(1-\frac{k^2}{n^2}\Big)^{1/2} < k^{1/2},
\ee
and
$$
  L_2 :=
  \Big(\frac{2n}{n+k}\Big)^{n+k} \Big(\frac{n-k}{2n}\Big)^{n-k}
 = 2^{2k} \frac{ (1 - \frac{k}{n})^{n-k}} {(1 + \frac{k}{n})^{n+k} }\,.
$$
We use then the inequalities
$(1-\frac{1}{x})^{x-1/2} <  \frac{1}{e}$ and
$(1+\frac{1}{x})^{x+1/2} > e$, where $x > 1$,
to derive
\baa
 &&    \Big(1 - \frac{k}{n}\Big)^{n-k}
 =  \Big(1 - \frac{k}{n}\Big)^{(\frac{n}{k}-\frac{1}{2})k}
   \Big(1 - \frac{k}{n}\Big)^{-k/2}
< \frac{1}{e^k} \Big(1 - \frac{k}{n}\Big)^{-k/2}\,, \\
&& \Big(1 + \frac{k}{n}\Big)^{n+k}
 = \Big(1 + \frac{k}{n}\Big)^{(\frac{n}{k}+\frac{1}{2})k}
   \Big(1 + \frac{k}{n}\Big)^{k/2}
 >\; e^k \Big(1 + \frac{k}{n}\Big)^{k/2}\,.
\eaa
Therefore
$$
   L_2
 < 2^{2k} \frac{1}{e^{2k}}
        \frac{1}{(1-\frac{k^2}{n^2})^{k/2}}\,,
$$
and that combined with \rf[es_k] proves \rf[tauk].

If $n \ge k^{3/2}$ and $k\geq 2$, then $(1- \frac{k^2}{n^2})^{k/4} >
(1 - \frac{1}{k})^{k/4} > 2^{-1/2}$, so \rf[tauk1] is valid with
$c_2 = 2^{1/2} c_1$. If $k=1$, then $n \ge 3$, and
$(1- \frac{k^2}{n^2})^{k/4} > 2^{-1/2}$ as well,
and that proves \rf[tauk1] as
well. \qed

\begin{lemma}
If $n-k = m$ is fixed and $n$ grows, then
\be
    \tau_{n,n-m}
\le c_3\, m^{1/4}
    \Big(\frac{me}{2}\Big)^{m/2} n^{-m/2}\,,\quad
     \label{k+m'}
\ee
\end{lemma}

\proof
We consiser the inequality \rf[p3]
$$
    \tau_{n,k}^2
\le c_1^2 k^{\frac{1}{2}}\,\Big(1-\frac{k^2}{n^2}\Big)^{1/2}
       \Big(\frac{2n}{n+k}\Big)^{n+k}
          \Big(\frac{n-k}{2n}\Big)^{n-k}\,,
$$
and then estimate the factors using substution $n-k = m$ where appropriate.
We have
\baa
   k^{\frac{1}{2}}\Big(1-\frac{k^2}{n^2}\Big)^{1/2}
&=& (n-k)^{1/2} \Big(\frac{k(n+k)}{n^2}\Big)^{1/2}
 \le 2^{1/2} m^{1/2}\,, \\
     \Big(\frac{2n}{n+k}\Big)^{n+k}
&=& \Big(1 + \frac{n-k}{n+k}\Big)^{n+k}
= \Big(1 + \frac{m}{n+k}\Big)^{n+k}
< e^m \,, \\
   \Big(\frac{n-k}{2n}\Big)^{n-k}
&=& \Big(\frac{m}{2n}\Big)^{m}\,.
\eaa
Thus, \rf[k+m'] follows with $c_3 = 2^{1/4} c_1$.
\qed

\begin{lemma}
If $k = \lfloor \la n \rfloor$,
where $\la \in (0,1)$, then as $n$ grows, we have an exponential decay
\be \lb{lan}
     \tau_{n,k}
\le c_4 n^{1/4} \rho_\la^{n/2}, \qquad \rho_\la < 1\,,
\ee
in particular
$$
   \tau_{n,n/2} < c_1 n^{1/4} \Big(\frac{4}{\sqrt{27}}\Big)^{n/2}.
$$
\end{lemma}

\proof
With $k = \lfloor \la n \rfloor$, set $\la' := \frac{k}{n}$, and note
that
\be \lb{la'}
     \la n - 1 \le k \le \la n
\ra  \la - \frac{1}{n} \le \la' \le \la\,.
\ee
Substitution $k = \la' n$ in \rf[p3] gives
\baa
      \tau_{n,k}^2
&\le& c_1^2\, k^{\frac{1}{2}}\,\Big(1-\frac{k^2}{n^2}\Big)^{1/2}
       \Big(\frac{2n}{n+k}\Big)^{n+k}
          \Big(\frac{n-k}{2n}\Big)^{n-k} \\
& < & c_1^2\, n^{1/2} \rho_{\la'}^n\,,
\eaa
where
$$
   \rho_{\la'}
 = \Big(\frac{2}{1+\la'}\Big)^{1+\la'}
   \Big(\frac{1-\la'}{2}\Big)^{1-\la'} < 1, \qquad \la' \in (0,1).
$$
On using that $g(x) := \ln\rho_x$ satisfies $g'(x) > -1$ for $x \in (0,1)$,
we derive from \rf[la'] that
$$
    \rho_{\la'} < e^{1/n} \rho_\la\,,
$$
and that proves \rf[lan] with $c_4 = e^{1/2n} c_1$.
If $\la = \frac{1}{2}$ we obtain
$
   \rho_{1/2}
 = 2 (\frac{1}{2})^{1/2}/(\frac{3}{2})^{3/2} = \frac{4}{\sqrt{27}}.
$
\qed


\section{The asymptotic formulas} \lb{as}


In this section, we derive the asympotic formulas \rf[tau*]-\rf[tau**]
of Theorem~\ref{t1.4}.

1) We start with the asymptotic formula for
$$
   \tau_k^* := \lim_{n\to \infty} \tau_{n,k}\,.
$$

For $\alpha, \beta > -1$, we denote by $\{P_m^{(\alpha,\beta)}\}$
the sequence of Jacobi polynomials which are orthogonal
with respect to the weight
$w_{\alpha,\beta}(x) = (1-x)^{\alpha} (1+x)^{\beta}$, with the
standard normalization
\be \lb{P}
   P_m^{(\alpha,\beta)}(1) = {m+\alpha \choose m}.
\ee
Note that derivatives of the Chebyshev polynomials are related to
Jacobi polynomials in the following way,
\be \lb{T=P}
   T_n^{(k)} = c_{n,k} P_m^{(\nu,\nu)}\,, \qquad
   m = n-k, \qquad \nu = k-\frac{1}{2}\,.
\ee

We will use the asymptotic property of Jacobi polynomials which is
described in terms of Bessel functions (see \cite[sect.\,8.1]{GS}),
namely the following equality from \cite{V}
\be \lb{e4.1}
       \lim_{m\to\infty} m^{-\alpha}P_m^{(\alpha,\beta)}(y_{m,r})
     = \left({j_{\alpha+1,r}\over2}\right)^{-\alpha}
       J_{\alpha}(j_{\alpha+1,r})\,,
\ee
where $\,y_{m,r}\,$ is the point of the $r$-th local extremum of
$\,P_m^{(\alpha,\beta)}$ counted in decreasing order and
$\,j_{\nu,r}\,$ is the $\,r$-th positive zero of the Bessel function
$\,J_{\nu}$.

\begin{lemma}\label{l4.1}
We have
\be \lb{tau*1}
    \tau_k^*
  = \Gamma(\nu+1)
    \Big({j_{\nu+1,1}\over2}\Big)^{-\nu}
    |J_\nu(j_{\nu+1,1})|\,, \qquad \nu  = k - \Frac{1}{2}\,.
\ee
\end{lemma}

\proof
By \rf[P]-\rf[T=P], since $\omega_{n,k} = y_{m,1}$, we have
$$
     \tau_k^*
  = \lim_{m\to\infty}
     \frac{|P_m^{(\nu,\nu)}(y_{m,1})|}
          {P_m^{(\nu,\nu)}(1)}
  =  \frac
     {\lim\limits_{m\to\infty} m^{-\nu} |P_m^{(\nu,\nu)}(y_{m,1})|}
     {\lim\limits_{m\to\infty} m^{-\nu} {m+\nu \choose m} }
  = \frac{L_1}{L_2}\,.
$$
By \rf[e4.1],
$$
   L_1
 = \Big({j_{\nu+1,1}\over2}\Big)^{-\nu}
        |J_\nu(j_{\nu+1,1})|,
$$
while for the denominator we use
$$
   {m+\nu \choose m}
 = \frac{\Gamma(m+\nu+1)}{\Gamma(m+1)\Gamma(\nu+1)}, \qquad
   \lim_{m\to\infty}
    m^{-\nu}\frac{\Gamma(m+\nu+1)}{\Gamma(m+1)} = 1\,,
$$
to obtain
$$
    L_2 = 1/\Gamma(\nu+1)\,
$$
and that proves the lemma.
\qed

\begin{lemma}[\cite{LW}]
The first positive zero $\,j_{\nu,1}\,$ of the Bessel function
$\,J_{\nu}$ obeys the following
asymptotic expansion
\be \lb{j}
    j_{\nu,1}
 =  \nu + a \nu^{1/3} + \OO(\nu^{-1/3}), \qquad
    a = -i_1/2^{1/3} = 1.8557...
\ee
where $i_1$ is the first negative zero of the Airy function $\Ai(x)\,$.
\end{lemma}

\begin{lemma}
We have
\be \lb{J}
   J_\nu(j_{\nu+1,1})
= - \Big(\frac{2}{\nu}\Big)^{2/3}\Ai'(i_1) + \OO(\nu^{-1}).
\ee
\end{lemma}

\proof
We will need the asymptotic behavior of $\,J_\nu(\nu x)$ for large
(fixed) $\,\nu\,$ and $\,x \ge 1$ (that is, around
the first positive zero $j_{\nu,1}$), which is given by the following formula
(see \cite[Chapter\,11]{Ol} or \cite{LW}),
\be \label{e4.4}
    J_\nu(\nu x)
  = {\phi(z)\over\nu^{1/3}} \left[\Ai(\nu^{2/3}z)
        \big(1 + \OO(\nu^{-2})\big)
      + \frac{\Ai'(\nu^{2/3}z)}{\nu^{4/3}}
         \big(B_0(z) + \OO(\nu^{-2})\big)\right]\,,
\end{equation}
where $0 < B_0(z) \le B_0(0)$ for $z \le 0$ and
$$
     z = - \left({3\over2}\sqrt{x^2-1}
               - {3\over2}\sec^{-1}(x)\right)^{2/3},\qquad
   \phi(z)= \left({4z\over 1-x^2}\right)^{1/4}\,, \qquad x \ge 1\,.
$$

Let $\,x=1+\delta\,$, where $\,\delta= \OO(\nu^{-2/3})$ and
$\,\delta>0$. Then
$$
    \sec^{-1}(x)
  = \arccos\Big({1\over 1+\delta}\Big)
  = \sqrt{2\delta}\,\Big(1-\frac{5\delta}{12}+\OO(\delta^2)\Big)\,,
$$
whence we obtain for $\,z\,$ and $\,\phi(z)\,$
$$
      z  = -2^{1/3}\delta\Big(1+\OO(\delta)\Big)\,, \qquad
 \phi(z) =  2^{1/3} + \OO(\delta)\,.
$$
Substitution of these quantities in \rf[e4.4] yields
\be \lb{J1}
    J_\nu(\nu(1+\delta))
 = \Big(\frac{2}{\nu}\Big)^{1/3}
    \Big( \Ai\big(-\nu^{2/3}2^{1/3}\delta\big)  + \OO(\nu^{-2/3})\Big)
   + \OO(\nu^{-1})
\ee
From \rf[j], we have
\baa
     j_{\nu+1,1}
 &=& \nu + 1 -\frac{i_1}{2^{1/3}} (\nu+1)^{1/3} + \OO(\nu^{-1/3}) \\
 &=&  \nu(1 + \delta_0), \qquad
      \delta_0
   = -\frac{i_1}{2^{1/3}} \nu^{-2/3} + \nu^{-1} + \OO(\nu^{-4/3})\,,
\eaa
so putting this into \rf[J1], we conclude
\baa
    J_\nu(j_{\nu+1,1})
&=& \Big(\frac{2}{\nu}\Big)^{1/3}
    \Big(\Ai\big(i_1 - (\Frac{2}{\nu})^{1/3} + \OO(\nu^{-2/3})\big)
         + \OO(\nu^{-2/3})\Big) \\
&=& \Big(\frac{2}{\nu}\Big)^{2/3}
    \Ai'(i_1) + \OO(\nu^{-1})\,,
\eaa
and that proves the lemma.
\qed

\medskip
{\bf Proof of Theorem \ref{t1.4}, part \rf[tau*].} With the
substitution $\nu = k - \frac{1}{2}$, we obtain \baa
    \big|J_{k-{1\over2}}(j_{k+{1\over2},1})\big|
&\stackrel{\rf[J]}{=}&
    \Big({2 \over k}\Big)^{2/3}\Big(|\Ai'(i_1)|+\OO(k^{-1/3})\Big)\,, \\
    j_{k+\frac{1}{2},1}
&\stackrel{\rf[j]}{=}&
     k +\Frac{1}{2} + ak^{1/3} + \OO(k^{-1/3}), \\
     \Big(\frac{j_{k+\frac{1}{2},1}}{2}\Big)^{-(k-\frac{1}{2})}
&=& \Big(\frac{2}{k}\Big)^{k-1/2} e^{-1/2} e ^{-ak^{1/3}}
       \big(1 + \OO(k^{-1/3})\big)\,, \\
     \Gamma(k+\Frac{1}{2})
& =& \frac{(2k)!}{4^k k!} \sqrt{\pi}
     \;=\; \Big(\frac{k}{e}\Big)^k \sqrt{2\pi} \big(1 + \OO(k^{-1})\big)\,,
\eaa
and formula \rf[tau*1] gives
$$
    \tau_k^*
 =  C_0 \Big({2\over e}\Big)^k e^{-a_0 k^{1/3}} k^{-1/6}
    \big(1+\OO(k^{-1/3})\big),
$$
where
$$
   C_0 = 4^{1/3}\sqrt{\pi\over e}\,|\Ai'(i_1)|\,,
   \qquad a_0 = - i_1/2^{1/3}\,,
$$
and that proves the first part of Theorem~\ref{t1.4}.
\qed

\medskip
2) Next, we will prove the asymptotic formula
for
$$
   \tau_m^{**} := \lim_{n\to \infty} n^{m/2} \tau_{n,n-m}\,.
$$
Note that, with $m = n-k$ fixed, and provided that the limit exists,
we have
$$
   \tau_m^{**}
:= \lim_{k\to \infty} (k+m)^{m/2} \tau_{k+m,k}
 = \lim_{k\to \infty} k^{m/2} \tau_{k+m,k}\,.
$$
We will use the relation
$$
    T_{k+m}^{(k)} = c_{m,k} P_m^{(\la)}, \qquad \la = k,
$$
and the asymptotic
properties of the ultraspherical polynomials $P_m^{(\la)}$
expressed in terms of the Hermite polynomials $\,H_m\,$  (see
\cite[eq.\,(5.6.3)]{GS})
\be \lb{H}
    \lim_{\lambda\to\infty} \lambda^{-m/2}
    P_m^{(\lambda)}\Big({x\over\sqrt\lambda}\Big)
  = {H_m(x)\over m!}\,.
\ee

\begin{lemma}
We have
$$
   \tau_m^{**} = 2^{-m}|H_m(x'_m)|\,,
$$
where $x'_m$ is the point of the rightmost extremum of $H_m$.
\end{lemma}

\proof
With $\omega_{k+m,k}$ and $x'_m$ being the points of
the rightmost local extrema of $T_{k+m}^{(k)} = c_{m,k} P_m^{(k)}$ and $H_m$,
respectively, it follows from \rf[H] that, for a fixed $m$, we have
$$
    \tau_{k+m,k}
  = {|P_m^{(k)}(\omega_{k+m,k})|\over P_m^{(k)}(1)}
\sim {k^{m/2} |H_m(x'_m)|\over m!{m+2k-1\choose m}}
\sim 2^{-m} k^{-m/2} |H_m(x'_m)|, \qquad k \to \infty,
$$
and this implies
$$
    \tau_m^{**}
  =  \lim_{k\to\infty}  k^{m/2}\tau_{k+m,k}
  = 2^{-m} |H_m(x'_m)|\,.
$$
Lemma is proved.
\qed

\medskip
{\bf Proof of Theorem \ref{t1.4}, part \rf[tau**].}
For approximation of $\,H_m(x_m')\,$ we will use the formula of
Plancherel - Rotach (\cite[Theorem\,8.22.9]{GS}). Actually, we need
only the third part of this theorem, concerning the approximation of
$\,H_m\,$ around its turning point, where the behaviour of the
polynomial changes from oscillatory to monotonically increasing.
It states that if
\be \lb{x}
   x = (2m+1)^{1\over2} - 2^{-{1\over2}}\,3^{-{1\over3}}\,m^{-{1\over6}}\,t,
      \qquad t\in \mathbb{C}, \\
\ee
then
\be \label{e4.7}
    e^{-x^2/2}\,H_m(x)
  = 3^{1\over3}\,\pi^{-{3\over4}}\,
      2^{{m\over2}+{1\over4}} (m!)^{1\over2}\,m^{-{1\over12}}
      \left\{A(t) + O(m^{-{2\over3}})\right\},
\ee
where $\,A(z)= 3^{-{1\over3}}\pi \Ai(-3^{-{1\over3}}z)\,$
is the normalized Airy
function. Moreover, the asymptotic formula \rf[e4.7] holds
uniformly when $\,t\in \mathbb{C}\,$ is bounded.

Let $\,x_m\,$ be the largest zero of $\,H_m\,$, then (\cite[eq.
(6.32.5)]{GS})
$$
    x_m
  = (2m+1)^{1\over2} -2^{-{1\over2}}\,3^{-{1\over3}}\,
    m^{-{1\over6}}\, i_1^* + \OO(m^{-5/6}), \qquad m\geq1,
$$
where $i_1^* = -3^{1/3}i_1\,$ is the first zero of $\,A(z)$. Since
$\,H'_m(x)=2m H_{m-1}(x)\,$, we have for $\,m\geq2\,$
\baa
     x_m'
 &=& x_{m-1} = (2m-1)^{1\over2}
     -2^{-{1\over2}}3^{-{1\over3}}m^{-{1\over6}} i^*_1 + \OO(m^{-5/6})\\
& =& x_m - (2m)^{-{1\over2}} + \OO(m^{-5/6}),
\eaa
and we can put $x_m'$ in the form \rf[x], with $t = t'_m$ where
$$
   t_m' = i^*_1 + 3^{1\over3} m^{-{1\over3}} + O(m^{-2/3}).
$$
Then formula \rf[e4.7] gives
\baa
      H_m(x_m')
  &=& e^{\frac{1}{2}(x_m')^2} 3^{1\over3}\,
      \pi^{-{3\over4}}\, 2^{{m\over2}+{1\over4}}\,
      (m!)^{1\over2}\, m^{-{1\over12}}
      \left\{A(t_m') + \OO(m^{-{2\over3}})\right\} \\
&=& -(2em)^{m\over2}\,e^{-|i_1|m^{1/3}}m^{-1/6}\sqrt{2\pi/e}\,
      \Ai'(i_1)\Big(1 + \OO(m^{-1/3})\Big).
\eaa
Finally, we obtain
$$
    \tau_m^{**} = 2^{-m}|H_m(x_m')|
    = \left(em\over2\right)^{m/2}
      e^{-a_1 m^{1/3}}m^{-1/6}\Big(1 + \OO(m^{-1/3})\Big).
$$
where
$$
   C_1 = \sqrt{2\pi\over e}\,\Ai'(i_1), \qquad a_1 = |i_1|\,.
$$
Theorem \ref{t1.4} is proved. \qed


\section{Remarks}


\textbf{1.} As was mentioned in introduction, Theorem~\ref{t2.1} is
due to Sz\'{a}sz \cite{os}. The statement of Theorem~\ref{t2.1}
appears in \cite[pp. 304--305]{aar} along with a proof of the
Legendre case ($\lambda=1/2$). The proof of this case originally was
given by Szeg\H{o}, who confirmed a conjecture made by J. Todd. Our
proof follows the same approach. An alternative proof of
Theorem~\ref{t1.1} can be obtained using some results of Bojanov and
Naidenov in \cite{BN}. \smallskip

\textbf{2.} To obtain better approximation of $\,\tau_{n,k}\,$ one
needs more precise asymptotic formulae for ultraspherical and
Hermite polynomials and bounds for their extreme zeros. In this
connection we refer to \cite{DJ, AE, EL, EL1, G, Kr, WZ}.
\bigskip\bigskip

\textbf{Acknowledgements.} The first two authors are supported by
the University of Sofia Research Fund under Contract 80-10-11/2017.
The second author acknowledges the support by the Bulgarian National
Research Fund through Contract DN 02/14. The third author is
supported by a Research Grant from Pembroke College, Cambridge.


\bigskip\noindent
{\sc Nikola Naidenov, Geno Nikolov} \smallskip\\
Department of Mathematics and Informatics\\
Sofia University "St. Kliment Ohridski\\
5 James Bourchier Blvd. \\
1164 Sofia \\
BULGARIA \smallskip \\
{\it E-mails:} $\begin{array}{l}
\text{nikola@fmi.uni-sofia.bg}\\
\text{geno@fmi.uni-sofia.bg}
        \end{array}$
\bigskip

\noindent
{\sc Alexei Shadrin} \smallskip\\
Department of Applied Mathematics\\
and Theoretical Physics (DAMTP) \\
Cambridge University \\
Wilberforce Road \\
Cambridge CB3 0WA \\
UNITED KINGDOM \\
{\it E-mail:} {\tt a.shadrin@damtp.cam.ac.uk}

\end{document}